\documentclass[10pt]{amsart}
\usepackage{amsmath,amssymb,latexsym,hyperref,graphicx}

\newtheorem*{corollary*}{Corollary}
\newtheorem*{proposition*}{Proposition}
\newtheorem*{theorem*}{Theorem}
\newtheorem{theorem}{Theorem}

\newcommand{\nequiv}{\mathrel{\not\equiv}}
\newcommand{\nonoverlappingsubwords}[2]{#2/#1}
\newcommand{\vphantp}{\vphantom{(}}
\newcommand{\vphantn}{\vphantom{\nonoverlappingsubwords{P}{n}}}

\begin{document}

\title{The number of nonzero \\ binomial coefficients modulo $p^\alpha$}
\author{Eric Rowland}

\address{
	Mathematics Department \\
	Tulane University \\
	New Orleans, LA 70118, USA
}
\date{March 19, 2011}

\thanks{I would like to thank Elizabeth Kupin and Doron Zeilberger for helpful discussions.}

\begin{abstract}
In 1947 Fine obtained an expression for the number $a_p(n)$ of binomial coefficients on row $n$ of Pascal's triangle that are nonzero modulo $p$.
In this paper we use Kummer's theorem to generalize Fine's theorem to prime powers, expressing the number $a_{p^\alpha}(n)$ of nonzero binomial coefficients modulo $p^\alpha$ as a sum over certain integer partitions.
For fixed $\alpha$, this expression can be rewritten to show explicit dependence on the number of occurrences of each subword in the base-$p$ representation of $n$.
\end{abstract}

\maketitle
\markboth{Eric Rowland}{The number of nonzero binomial coefficients modulo $p^\alpha$}

\section{Introduction}\label{Introduction}

The study of arithmetic properties of binomial coefficients has a rich history.
A main theme is that properties of $\binom{n}{m}$ modulo a prime $p$ are related to the base-$p$ representations of $n$ and $m$.
Let $n_l n_{l-1} \cdots n_0$ be the word consisting of the standard base-$b$ digits of a nonnegative integer $n$.
We use $n$ and $n_l n_{l-1} \cdots n_0$ interchangeably.
We consider the base-$b$ representation of $0$ to be the empty word $\epsilon$.
For $0 \leq m \leq n$ we write $m_l m_{l-1} \cdots m_0$ for the base-$b$ representation of $m$, where we pad with zeros if it is otherwise shorter than $n_l n_{l-1} \cdots n_0$.
With the exception of the proposition in Section~\ref{Generalizing}, we will take $b = p$ to be prime.

Two of the classic results are Kummer's theorem of 1852~\cite[pages~115--116]{Kummer} and Lucas' theorem of 1878~\cite{Lucas}.

\begin{theorem*}[Kummer]
Let $p$ be a prime, and let $0 \leq m \leq n$.
The exponent of the highest power of $p$ dividing $\binom{n}{m}$ is the number of borrows involved in subtracting $m$ from $n$ in base $p$.
\end{theorem*}

\begin{theorem*}[Lucas]
Let $p$ be a prime, and let $0 \leq m \leq n$.
Then
\[
	\binom{n}{m} \equiv \prod_{i=0}^l \binom{n_i}{m_i} \mod p.
\]
\end{theorem*}

Let $a_k(n)$ be the number of integers $0 \leq m \leq n$ such that $\binom{n}{m} \nequiv 0 \mod k$; that is, $a_k(n)$ is the number of nonzero entries on row $n$ of Pascal's triangle modulo $k$.
Let $|n|_w$ be the number of occurrences of the word $w$ in $n_l n_{l-1} \cdots n_0$.

In 1899 Glaisher~\cite[\S14]{Glaisher} initiated the study of counting entries on row $n$ of Pascal's triangle modulo $k$ by using Lucas' theorem to determine $a_2(n) = 2^{|n|_1}$.
The proof is simple: In order that $\binom{n}{m}$ be odd, each term $\binom{n_i}{m_i}$ in the product must be $1$, so if $n_i = 0$ then $m_i = 0$ and if $n_i = 1$ then $m_i$ can be either $0$ or $1$.

In 1947 Fine~\cite{Fine} generalized Glaisher's result to an arbitrary prime.
Fine's result follows from Lucas' theorem in the same way.

\begin{theorem*}[Fine]
Let $p$ be a prime, and let $n \geq 0$.
The number of nonzero entries on row $n$ of Pascal's triangle modulo $p$ is
\[
	a_p(n) = \prod_{i=0}^l \left(n_i + 1\right).
\]
\end{theorem*}

Note that Fine's expression may be rewritten as
\[
	a_p(n) = \prod_{r=0}^{p-1} (r+1)^{|n|_r},
\]
which more directly identifies the contribution of each digit $0 \leq r \leq p-1$.

In Section~\ref{Generalizing} we generalize Fine's result to prime powers, obtaining a formula for $a_{p^\alpha}(n)$.
In Section~\ref{Expressions} we provide an algorithm for rewriting this formula in terms of $|n|_w$, as we have just done with Fine's expression.
Previously, such formulas were only known for $a_4(n)$, $a_9(n)$, and $a_8(n)$.

We mention that one may generalize Glaisher's result in a different direction, namely to ask for the number $a_{k,r}(n)$ of integers $0 \leq m \leq n$ such that $\binom{n}{m} \equiv r \mod k$.
In this context, Fine's result is an evaluation of the sum over all nonzero residue classes when $k = p$ is prime, and the main result of this paper is an evaluation of the sum
\[
	a_{p^\alpha}(n) = \sum_{r=1}^{p^\alpha - 1} a_{p^\alpha,r}(n)
\]
for a prime power modulus.

There have been several studies of $a_{k,r}(n)$.
For prime $k = p$, Hexel and Sachs~\cite[\S5]{Hexel--Sachs} determined a formula for $a_{p,r^i}(n)$ in terms of $(p-1)$th roots of unity, where $r$ is a primitive root modulo $p$, and from this obtained $a_{3,1}(n) = 2^{|n|_1 - 1} (3^{|n|_2} + 1)$, $a_{3,2}(n) = 2^{|n|_1 - 1} (3^{|n|_2} - 1)$, and explicit formulas for $a_{5,r^i}(n)$ in terms of $|n|_1$, $|n|_2$, $|n|_3$, and $|n|_4$.
Garfield and Wilf~\cite{Garfield--Wilf} provided an algorithm to compute the generating function $\sum_{i=0}^{p-2} a_{p,r^i}(n) x^i$, where again $r$ is a primitive root.
Recently, Amdeberhan and Stanley~\cite[Theorem~2.1]{Amdeberhan--Stanley} studied the number of coefficients equal to $r$ in the $n$th power of a general multivariate polynomial over a finite field, where $r$ is an invertible element of the field.

In the late 1980s researchers began to consider $a_{k,r}(n)$ for certain prime power moduli $k = p^\alpha$.
Davis and Webb~\cite{Davis--Webb 1989} gave formulas for $a_{4,1}(n)$, $a_{4,2}(n)$, and $a_{4,3}(n)$ in terms of $|n|_1$, $|n|_{10}$, and $|n|_{11}$.
Around the same time, Granville~\cite{Granville 1992} showed that if $r$ is odd and $\alpha \in \{2, 3\}$ then $a_{2^\alpha,r}(n)$ is either $0$ or a power of $2$.
Huard, Spearman, and Williams~\cite{Huard--Spearman--Williams mod 9, Huard--Spearman--Williams mod 8} gave formulas for $a_{9,r}(n)$ and $a_{8,r}(n)$ (and their sums, which we derive again below).
Some of these results use a generalization of Lucas' theorem to prime powers found by Davis and Webb~\cite{Davis--Webb 1990}.

There has also been some general work on squares of primes.
Huard, Spearman, and Williams~\cite{Huard--Spearman--Williams mod p^2} used the result of Hexel and Sachs to find, when $p \mid r$ and $r \neq 0$, a formula for $a_{p^2,r}(n)$ depending only on $|n|_w$ for words $w$ of length at most $2$.
Earlier, Webb~\cite[Theorem~3]{Webb} showed if $p \nmid r$ then $a_{p^2,r}(n)$ does \emph{not} depend only the subwords of length at most $2$ but does depend only on the blocks of nonzero digits in $n$.
However, the corollary in the next section implies that by summing $a_{p^2,r}(n)$ over all nonzero residue classes $r$ modulo $p^2$ the dependence on only the subwords of length at most $2$ is achieved.

We would be remiss to not mention Granville's thorough survey~\cite{Granville 1997}, which discusses many additional arithmetic aspects of binomial coefficients and provides another generalization of Lucas' theorem to prime powers.

\section{Generalizing Fine's theorem}\label{Generalizing}

We adopt the usual conventions that an empty sum is $0$, an empty product is $1$, and there is precisely one integer partition of $0$ (namely, the empty set).

For $k \geq 1$, let
\[
	c(w_k w_{k-1} \cdots w_0) = \frac{w_k}{w_k + 1} \cdot \left(\prod_{h=1}^{k-1} \frac{b - w_{k-h}}{w_{k-h} + 1} \right) \cdot \frac{b - w_0 - 1}{w_0 + 1}.
\]
The function $c$ assigns a rational number to a word on the alphabet $\{0, 1, \dots, b-1\}$.
We will see this function arise naturally in the proof of Theorem~\ref{main}.

For a nonnegative integer $\gamma$, let $S_\alpha(\gamma)$ be the set of integer partitions of $\gamma$ into at least $\max(0, \gamma - (\alpha - 1))$ parts, all of size at least $2$.
For example,
\[
	S_8(10) = \{\{6,2,2\}, \{5,3,2\}, \{4,4,2\}, \{4,3,3\}, \{4,2,2,2\}, \{3,3,2,2\}, \{2,2,2,2,2\}\}.
\]

Let $|P|$ be the number of parts in the integer partition $P$.
Let $$\sum_{\nonoverlappingsubwords{P}{n}} c(v) c(w) \cdots c(z)$$ be the sum over the sets $\{v, w, \dots, z\}$ of $|P|$ nonoverlapping subwords of $n = n_l n_{l-1} \cdots n_0$ such that the multiset $\{|v|, |w|, \dots, |z|\}$ of subword lengths is equal to $P$.
In this sum we consider two subwords $n_{i_1} n_{i_1-1} \cdots n_{f_1}$ and $n_{i_2} n_{i_2-1} \cdots n_{f_2}$ to be distinct precisely when $i_1 \neq i_2$ or $f_1 \neq f_2$, so it would be more precise (but more cumbersome) to say that $\sum_{\nonoverlappingsubwords{P}{n}}$ is a sum over certain sets of pairs of indices.

For example, if $n = n_5 n_4 n_3 n_2 n_1 n_0$ then
\begin{align*}
	\sum_{\nonoverlappingsubwords{\{3,2\}}{n}} &c(v) c(w) \cdots c(z) \\
	&= c(n_5 n_4 n_3) c(n_2 n_1) + c(n_5 n_4 n_3) c(n_1 n_0) + c(n_4 n_3 n_2) c(n_1 n_0) \\
	&+ c(n_5 n_4) c(n_3 n_2 n_1) + c(n_5 n_4) c(n_2 n_1 n_0) + c(n_4 n_3) c(n_2 n_1 n_0)
\end{align*}
and
\begin{align*}
	\sum_{\nonoverlappingsubwords{\{2,2\}}{n}} &c(v) c(w) \cdots c(z) \\
	&= c(n_5 n_4) c(n_3 n_2) + c(n_5 n_4) c(n_2 n_1) + c(n_5 n_4) c(n_1 n_0) \\
	&+ c(n_4 n_3) c(n_2 n_1) + c(n_4 n_3) c(n_1 n_0) + c(n_3 n_2) c(n_1 n_0).
\end{align*}

We now have the notation to state the main result of the paper.

\begin{theorem}\label{main}
Let $p$ be a prime, let $\alpha \geq 0$, and let $n \geq 0$.
The number of nonzero entries on row $n$ of Pascal's triangle modulo $p^\alpha$ is
\[
	a_{p^\alpha}(n) = \left(\prod_{\vphantn i=0}^l \left(n_i + 1\right)\right) \sum_{\vphantn \gamma = 0}^{2 (\alpha - 1)} \sum_{P \in S_\alpha(\gamma)} \sum_{\nonoverlappingsubwords{P}{n}} c(v) c(w) \cdots c(z).
\]
\end{theorem}

Note that if it is convenient we may extend the sum over $S_\alpha(\gamma)$ to a sum over partitions including $1$ if we set $c(w_0) = 0$.

Theorem~\ref{main} follows from the following proposition.
Let $b \geq 2$, and let $A_n(\beta)$ be the number of integers $0 \leq m \leq n$ such that there are exactly $\beta$ borrows involved in computing $n - m$ in base $b$.
Let $S(\gamma, \delta)$ be the set of integer partitions of $\gamma$ into $\delta$ parts where all parts are at least $2$.
For example, $S(10, 3) = \{\{6,2,2\}, \{5,3,2\}, \{4,4,2\}, \{4,3,3\}\}$.

\begin{proposition*}
Let $b \geq 2$, let $\alpha \geq 0$, and let $n \geq 0$.
Then
\[
	\frac{A_n(\beta)}{A_n(0)} = \sum_{\vphantn \gamma = \beta}^{2 \beta} \sum_{\vphantn P \in S(\gamma, \gamma-\beta)} \sum_{\nonoverlappingsubwords{P}{n}} c(v) c(w) \cdots c(z).
\]
\end{proposition*}

Everett~\cite{Everett 2008} gave a different expression for $A_n(\beta)$ as a sum over all length-$l$ words on $\{0, 1\}$ with precisely $\beta$ $1$s.
Everett's expression is simpler to state and faster to compute for an explicit integer $n$.
However, because of its high-level dependence on $l$, it is farther away from being able to produce formulas in terms of subword counts.
Note that neither expression for $A_n(\beta)$ relies on the base $b$ being prime.

Now let $b = p$ be prime.
By Kummer's theorem, $\binom{n}{m} \nequiv 0 \mod p^\alpha$ precisely when there are fewer than $\alpha$ borrows when subtracting $m$ from $n$ in base $p$.
Therefore $a_{p^\alpha}(n) = \sum_{\beta=0}^{\alpha-1} A_n(\beta)$.
Substituting the expression for $A_n(\beta)/A_n(0)$ in the proposition and interchanging the two outermost sums gives the statement of the theorem.

Therefore it suffices to prove the proposition.
For $n = n_l n_{l-1} \cdots n_0$ and $m = m_l m_{l-1} \cdots m_0$, let $n' = n_{l-1} \cdots n_0$ and $m' = m_{l-1} \cdots m_0$.
Furthermore, let $n^{(i)} = n_{l-i} \cdots n_0$.

\begin{proof}[Proof of the proposition]
We first find a recurrence for $A_n(\beta)$ by establishing the relationship between borrows in $n - m$ and borrows in $n' - m'$.
Since $A_n(\beta) = 0$ when $\beta > l$, it suffices to consider $\beta \leq l$.

It may happen that $m' > n'$ even if $m \leq n$, so we must decide how to count borrows in the computation of $n' - m'$ in this case.
The standard subtraction algorithm produces infinitely many borrows.
However, the only borrows that are preserved when passing from $n' - m'$ to $n - m$ are those up through the borrow from the $l$th digit in $n'$ (which is $0$).
Therefore, let $B_n(\beta)$ be the number of integers $n < m \leq b^{l+1} - 1$ such that there are exactly $\beta$ borrows up through the borrow from $n_{l+1} = 0$ involved in computing $n - m$.

Now we write $A_n(\beta)$ in terms of $A_{n'}(\beta)$ and $B_{n'}(\beta)$.
In the computation of $n - m$, a borrow from the digit $n_{i+1}$ occurs if $m_i > n_i$.
Moreover, if there is a borrow from $n_i$ then the borrow is propagated to $n_{i+1}$ whenever $m_i > n_i - 1$.
Thus if $m' \leq n'$ then there are $n_l + 1$ choices for $m_l$ such that $m \leq n$.
Similarly, if $m' > n'$ then there are $n_l$ choices for $m_l$ such that $m \leq n$.
Therefore
\[
	A_n(\beta) = (n_l + 1) A_{n'}(\beta) + n_l B_{n'}(\beta).
\]

We find a recurrence for $B_n(\beta)$ analogously:
If $m' \leq n'$ then there are $b - n_l - 1$ choices for $m_l$ such that $m > n$.
If $m' > n'$ then there are $b - n_l$ choices for $m_l$ such that $m > n$.
In each case we gain one additional borrow, so
\[
	B_n(\beta) = (b - n_l - 1) A_{n'}(\beta-1) + (b - n_l) B_{n'}(\beta-1).
\]

Iteratively substituting the equation for $B_n(\beta)$ into the equation for $A_n(\beta)$ until we reach $B_{n^{(\beta+1)}}(0) = 0$ produces the recurrence
\[
	A_n(\beta) = (n_l + 1) A_{n'}(\beta) + \sum_{i=1}^\beta n_l \left(\prod_{j=1}^{i - 1} (b - n_{l - j}) \right) (b - n_{l - i} - 1) A_{n^{(i+1)}}(\beta - i)
\]
for $n \geq 1$ and $0 \leq \beta \leq l$.
Divide both sides of this recurrence by $A_n(0)$, which (as in Fine's theorem) is $\prod_{i=0}^l \left(n_i + 1\right)$, to obtain
\[
	\frac{A_n(\beta)}{A_n(0)} - \frac{A_{n'}(\beta)}{A_{n'}(0)} = \sum_{i=1}^\beta c(n_l n_{l-1} \cdots n_{l-i}) \cdot \frac{A_{n^{(i+1)}}(\beta - i)}{A_{n^{(i+1)}}(0)},
\]
where $c(w_k w_{k-1} \cdots w_0)$ is as defined above.
Replacing $n$ with $n^{(j)}$ in this equation and summing over $0 \leq j \leq l - \beta$ causes the left side to telescope, and we see that
\[
	\frac{A_n(\beta)}{A_n(0)} - \frac{A_{n^{(l-\beta+1)}}(\beta)}{A_{n^{(l-\beta+1)}}(0)} = \sum_{j=0}^{l-\beta} \sum_{i=1}^\beta c(n_{l-j} n_{l-j-1} \cdots n_{l-j-i}) \cdot \frac{A_{n^{(j+i+1)}}(\beta-i)}{A_{n^{(j+i+1)}}(0)}.
\]
If $\beta \geq 1$ then $A_{n^{(l-\beta+1)}}(\beta) = 0$, and if $\beta = 0$ then $A_{n^{(l-\beta+1)}}(\beta) = A_\epsilon(0) = 1$.

We now verify that the expression for $A_n(\beta)/A_n(0)$ given in the statement of the proposition satisfies this recurrence and the correct boundary conditions.
The boundary conditions are easily checked; for $\beta = 0$ the expression is $1$, and for $\beta > l$ it is $0$.
After substituting, the right side of the recurrence is
\begin{align*}
	\sum_{\vphantn j=0}^{\vphantp l-\beta} &\sum_{\vphantn i=1}^{\vphantp \beta} c(n_{l-j} \cdots n_{l-j-i}) \sum_{\vphantn \gamma = \beta - i}^{2 (\beta - i)} \sum_{\vphantn P \in S(\gamma, \gamma-\beta+i)} \sum_{\nonoverlappingsubwords{P}{n^{(j+i+1)}}} c(v) c(w) \cdots c(z) \\
	&= \sum_{\vphantn j=0}^{l-\beta} \sum_{\vphantn i=1}^\beta \sum_{\vphantn \gamma = \beta + 1}^{2 \beta + 1 - i} \sum_{\vphantn P \in S(\gamma-1-i, \gamma-1-\beta)} \sum_{\nonoverlappingsubwords{P}{n^{(j+i+1)}}} c(n_{l-j} \cdots n_{l-j-i}) c(v) c(w) \cdots c(z) \\
	&= \sum_{\vphantn \gamma = \beta + 1}^{2 \beta} \sum_{\vphantn i=1}^{2 \beta + 1 - \gamma} \sum_{\vphantn j=0}^{l-\beta} \sum_{\vphantn P \in S(\gamma-1-i, \gamma-1-\beta)} \sum_{\nonoverlappingsubwords{P}{n^{(j+i+1)}}} c(n_{l-j} \cdots n_{l-j-i}) c(v) c(w) \cdots c(z)
\end{align*}
after shifting $\gamma \mapsto \gamma - 1 - i$ and interchanging the sums over $i$ and $\gamma$.

Momentarily fix $\beta + 1 \leq \gamma \leq 2 \beta$.
For each $1 \leq i \leq 2 \beta + 1 - \gamma$, take each integer composition of $\gamma - 1 - i$ into $\gamma - 1 - \beta$ parts, where all parts are at least $2$, and prepend $i+1$ to get a composition of $\gamma$ into $\gamma - \beta$ parts at least $2$.
In doing this we form each composition of $\gamma$ into $\gamma - \beta$ parts at least $2$ precisely once.
Therefore
\begin{align*}
	\sum_{\vphantn i=1}^{2 \beta + 1 - \gamma} &\sum_{\vphantn j=0}^{l-\beta} \sum_{\vphantn P \in S(\gamma-1-i, \gamma-1-\beta)} \sum_{\nonoverlappingsubwords{P}{n^{(j+i+1)}}} c(n_{l-j} \cdots n_{l-j-i}) c(v) c(w) \cdots c(z) \\
	&= \sum_{\vphantn P \in S(\gamma, \gamma-\beta)} \sum_{\nonoverlappingsubwords{P}{n}} c(v) c(w) \cdots c(z).
\end{align*}
The right side of the recurrence then becomes
\[
	\sum_{\vphantn \gamma = \beta + 1}^{2 \beta} \sum_{\vphantn P \in S(\gamma, \gamma-\beta)} \sum_{\vphantn \nonoverlappingsubwords{P}{n}} c(v) c(w) \cdots c(z),
\]
which if $\beta \geq 1$ is equal to
\[
	\sum_{\vphantn \gamma = \beta}^{2 \beta} \sum_{\vphantn P \in S(\gamma, \gamma-\beta)} \sum_{\nonoverlappingsubwords{P}{n}} c(v) c(w) \cdots c(z)
\]
and if $\beta = 0$ is equal to
\[
	\sum_{\vphantn \gamma = \beta}^{2 \beta} \sum_{\vphantn P \in S(\gamma, \gamma-\beta)} \sum_{\nonoverlappingsubwords{P}{n}} c(v) c(w) \cdots c(z) - 1 = 0
\]
as desired.
\end{proof}

We mention that the recurrences appearing early in the proof are sufficient to compute $a_{p^\alpha}(n)$ symbolically for fixed $\alpha$; for example, for $\beta = 0$ we have
\[
	A_n(0) = (n_l + 1) A_{n'}(0),
\]
giving Fine's theorem
\[
	a_p(n) = A_n(0) = \prod_{i=0}^l \left(n_i + 1\right).
\]
Of course, Fine's theorem also follows from the full statement of Theorem~\ref{main}; $S_1(0) = \{\{\}\}$, so the inner sum is a sum over one term, and the summand is the empty product.

For $\alpha = 2$ we have $S_2(0) = \{\{\}\}$, $S_2(1) = \{\}$, and $S_2(2) = \{\{2\}\}$, so
\[
	a_{p^2}(n) = \left(\prod_{i=0}^l \left(n_i + 1\right)\right) \cdot \left(\sum_{\nonoverlappingsubwords{\{\}}{n}} 1 + \sum_{\nonoverlappingsubwords{\{2\}}{n}} c(w)\right).
\]
The second sum is simply the sum over all subwords of length $2$, so we have proved the following corollary.

\begin{corollary*}
Let $p$ be a prime, and let $n \geq 0$.
The number of nonzero entries on row $n$ of Pascal's triangle modulo $p^2$ is
\[
	a_{p^2}(n) = \left(\prod_{i=0}^l \left(n_i + 1\right)\right) \cdot
		\left(1 + \sum_{i=0}^{l-1} \frac{n_{i+1}}{n_{i+1} + 1} \cdot \frac{p - n_i - 1}{n_i + 1} \right).
\]
\end{corollary*}

Brief words are in order regarding how one can experimentally guess the general expression for $A_n(\beta)/A_n(0)$ given by the proposition once one knows the recurrence
\[
	\frac{A_n(\beta)}{A_n(0)} - \frac{A_{n^{(l-\beta+1)}}(\beta)}{A_{n^{(l-\beta+1)}}(0)} = \sum_{j=0}^{l-\beta} \sum_{i=1}^\beta c(n_{l-j} n_{l-j-1} \cdots n_{l-j-i}) \cdot \frac{A_{n^{(j+i+1)}}(\beta-i)}{A_{n^{(j+i+1)}}(0)}.
\]
It is clear from this recurrence that the fully resolved expression for $A_n(\beta)/A_n(0)$ is a sum of terms of the form $c(v) c(w) \cdots c(z)$, where $(v, w, \dots, z)$ is a tuple of nonoverlapping subwords of $n$.
For example, if $n = n_5 n_4 n_3 n_2 n_1 n_0$ then
\begin{multline*}
	A_n(3)/A_n(0) = c(n_5 n_4 n_3 n_2) + c(n_4 n_3 n_2 n_1) + c(n_3 n_2 n_1 n_0) \\
	+ c(n_5 n_4 n_3) c(n_2 n_1) + c(n_5 n_4 n_3) c(n_1 n_0) + c(n_4 n_3 n_2) c(n_1 n_0) \\
	+ c(n_5 n_4) c(n_3 n_2 n_1) + c(n_5 n_4) c(n_2 n_1 n_0) + c(n_4 n_3) c(n_2 n_1 n_0) \\
	+ c(n_5 n_4) c(n_3 n_2) c(n_1 n_0).
\end{multline*}
Moreover, each tuple appears at most once.
Thus it suffices to determine which tuples appear.

Upon explicitly computing $A_n(3)/A_n(0)$ and several additional values, one observes that if a tuple $(v, w, \dots, z)$ appears in $A_n(\beta)/A_n(0)$ and $(\tilde{v}, \tilde{w}, \dots, \tilde{z})$ is a tuple of nonoverlapping subwords of the same length such that the multisets $\{|v|, |w|, \dots, |z|\}$ and $\{|\tilde{v}|, |\tilde{w}|, \dots, |\tilde{z}|\}$ are equal, then $(\tilde{v}, \tilde{w}, \dots, \tilde{z})$ also seems to appear, regardless of the order that the subwords of either tuple occur in $n$.
For example, all pairs of nonoverlapping subwords with lengths $\{3, 2\}$ appear in $A_n(3)/A_n(0)$.
So presumably it suffices to determine which multisets of subword lengths appear for a given $\beta$ (and $l$).
From the data, one guesses that the multisets are certain integer partitions of integers $\gamma$ into $\gamma - \beta$ parts, hence the proposition.
For example, the set of partitions appearing in $A_n(3)/A_n(0)$ is $\{\{4\}, \{3, 2\}, \{2, 2, 2\}\}$.

\section{Expressions in terms of $|n|_w$}\label{Expressions}

In this section we describe, for fixed $\alpha$, how to rewrite the expression for $a_{p^\alpha}(n)$ of the previous section to show explicit dependence on the subword counts for symbolic $n$.
Note that if $w$ begins with $0$ or ends with $p-1$ then $c(w) = 0$, so $a_{p^\alpha}(n)$ does not depend on $|n|_w$.

In Section~\ref{Introduction} we did this for Fine's theorem (where $\alpha = 1$).
For $\alpha = 2$ it is also done easily; collecting identical terms of the sum appearing in the corollary yields
\[
	a_{p^2}(n) = \left(\prod_{w_0 = 0}^{p-1} (w_0 + 1)^{|n|_{w_0}}\right) \cdot
		\left(1 + \sum_{w_1 = 0}^{p-1} \sum_{w_0 = 0}^{p-1} \frac{w_1}{w_1 + 1} \cdot \frac{p - w_0 - 1}{w_0 + 1} \cdot |n|_{w_1 w_0} \right).
\]
Now any explicit prime can be substituted to produce a formula.
For example, $p=2$ gives $a_4(n) = 2^{|n|_1} (1 + \frac{1}{2} |n|_{10})$.
For $p=3$ we have
\[
	a_9(n) = 2^{|n|_1} 3^{|n|_2} \left(1 + |n|_{10} + \frac{1}{4} |n|_{11} + \frac{4}{3} |n|_{20} + \frac{1}{3} |n|_{21}\right)
\]
(first found by Huard, Spearman, and Williams~\cite{Huard--Spearman--Williams mod 9}),
for $p=5$ we have
\begin{multline*}
	\frac{a_{25}(n)}{2^{|n|_1} 3^{|n|_2} 4^{|n|_3} 5^{|n|_4}} =
		1
		+ 2 |n|_{10}
		+ \frac{3}{4} |n|_{11}
		+ \frac{1}{3} |n|_{12}
		+ \frac{1}{8} |n|_{13} \\
		+ \frac{8}{3} |n|_{20}
		+ |n|_{21}
		+ \frac{4}{9} |n|_{22}
		+ \frac{1}{6} |n|_{23}
		+ 3 |n|_{30}
		+ \frac{9}{8} |n|_{31}
		+ \frac{1}{2} |n|_{32}
		+ \frac{3}{16} |n|_{33} \\
		+ \frac{16}{5} |n|_{40}
		+ \frac{6}{5} |n|_{41}
		+ \frac{8}{15} |n|_{42}
		+ \frac{1}{5} |n|_{43},
\end{multline*}
and so on.
To give a very explicit example, the base-$5$ representation of $1947$ is $30242$, so $a_{25}(1947) = 3^2 \cdot 4^1 \cdot 5^1 \cdot (1 + 3 + 8/15) = 816$.

For $\alpha = 3$ the theorem provides
\[
	\frac{a_{p^3}(n)}{a_p(n)} = \sum_{\nonoverlappingsubwords{\{\}}{n}} 1 + \sum_{\nonoverlappingsubwords{\{2\}}{n}} c(w) + \sum_{\nonoverlappingsubwords{\{3\}}{n}} c(w) + \sum_{\nonoverlappingsubwords{\{2,2\}}{n}} c(v) c(w).
\]
The first three sums can be directly rewritten in terms of $|n|_w$.
The final sum over nonoverlapping pairs of length-$2$ subwords can be written as the sum over unrestricted pairs of subwords minus the sum over overlapping pairs of subwords (of which there are two kinds --- overlapping in one letter and overlapping in both):
\begin{multline*}
	2 \sum_{\nonoverlappingsubwords{\{2,2\}}{n}} c(v) c(w) \\
	= \sum_{i=0}^{l-1} \sum_{j=0}^{l-1} c(n_{i+1} n_i) c(n_{j+1} n_j)
	- 2 \sum_{i=0}^{l-2} c(n_{i+2} n_{i+1}) c(n_{i+1} n_i)
	- \sum_{i=0}^{l-1} c(n_{i+1} n_i)^2.
\end{multline*}
The coefficients take care of symmetries among the subword lengths.
Since each of these three new sums consists of sums over the entire word length, they can be rewritten to show the dependence on subwords of lengths $2$ and $3$ as
\[
	\sum_{v \in [p]^2} \sum_{w \in [p]^2} c(v) c(w) |n|_v |n|_w
	- 2 \sum_{w \in [p]^3} c(w_2 w_1) c(w_1 w_0) |n|_w
	- \sum_{w \in [p]^2} c(w)^2 |n|_w,
\]
where $[p] = \{0, 1, \dots, p-1\}$ is the alphabet of base-$p$ digits.
Thus we can write out an expression for $a_{p^3}(n)$ in terms of $|n|_w$.
Letting $p=2$ in this expression gives
\[
	a_8(n) = 2^{|n|_1} \left(1 + \frac{3}{8} |n|_{10} + |n|_{100} + \frac{1}{4} |n|_{110} + \frac{1}{8} |n|_{10}^2\right),
\]
which was obtained by Huard, Spearman, and Williams~\cite{Huard--Spearman--Williams mod 8}.
Formulas for other primes can be found similarly: $p=3$ gives
\begin{multline*}
	\frac{a_{27}(n)}{2^{|n|_1} 3^{|n|_2}} =
		1
		+ \frac{1}{2} |n|_{10}
		+ \frac{7}{32} |n|_{11}
		+ \frac{4}{9} |n|_{20}
		+ \frac{5}{18} |n|_{21} \\
		+ 3 |n|_{100}
		+ \frac{3}{4} |n|_{101}
		+ \frac{3}{4} |n|_{110}
		+ \frac{3}{16} |n|_{111}
		+ \frac{1}{3} |n|_{120}
		+ \frac{1}{12} |n|_{121} \\
		+ 4 |n|_{200}
		+ |n|_{201}
		+ |n|_{210}
		+ \frac{1}{4} |n|_{211}
		+ \frac{4}{9} |n|_{220}
		+ \frac{1}{9} |n|_{221} \\
		+ \frac{1}{2} |n|_{10}^2
		+ \frac{1}{4} |n|_{10} |n|_{11}
		+ \frac{4}{3} |n|_{10} |n|_{20}
		+ \frac{1}{3} |n|_{10} |n|_{21}
		+ \frac{1}{32} |n|_{11}^2 \\
		+ \frac{1}{3} |n|_{11} |n|_{20}
		+ \frac{1}{12} |n|_{11} |n|_{21}
		+ \frac{8}{9} |n|_{20}^2
		+ \frac{4}{9} |n|_{20} |n|_{21}
		+ \frac{1}{18} |n|_{21}^2.
\end{multline*}

For a general $\alpha$ we will need to be able to rewrite $\sum_{\nonoverlappingsubwords{P}{n}} c(v) c(w) \cdots c(z)$ --- the sum over subwords with $0$ overlaps --- in terms of $|n|_w$ for any given partition $P$.
To do this we can use inclusion--exclusion to express this sum as sums over sets of words with forced overlap conditions rather than restrictive overlap conditions:
\[
	\sum_{\nonoverlappingsubwords{P}{n}} c(v) c(w) \cdots c(z)
	= \sum_{i \geq 0} (-1)^i \sum_{\substack{\text{ways to guarantee} \\ \text{$i$ overlaps}}} c(v) c(w) \cdots c(z).
\]
Now each term is a sum over all sets of subwords of the desired lengths where certain pairs of subwords are required to overlap.
For such a sum, determine the ``connected components'' induced by these pairs, and for each connected component find all clusters of the subwords in which the required pairs overlap.
Then allow each connected component to range independently over the entire word $n$.
Since each sum now is over the entire word $n$, the expression can readily be rewritten in terms of $|n|_w$.

An implementation of this procedure is available in the \textit{Mathematica} package \textsc{BinomialCoefficients}~\cite{BinomialCoefficients}.
This implementation produces formulas for $a_{16}(n)$, $a_{32}(n)$, and $a_{64}(n)$ fairly quickly on a standard machine (in less than a minute for $a_{64}(n)$).
For example,
\begin{multline*}
	\frac{a_{16}(n)}{2^{|n|_1}} =
		1
		+ \frac{5}{12} |n|_{10}
		+ \frac{1}{2} |n|_{100}
		+ \frac{1}{8} |n|_{110}
		+ 2 |n|_{1000}
		+ \frac{1}{2} |n|_{1010}
		+ \frac{1}{2} |n|_{1100} \\
		+ \frac{1}{8} |n|_{1110}
		+ \frac{1}{16} |n|_{10}^2
		+ \frac{1}{2} |n|_{10} |n|_{100}
		+ \frac{1}{8} |n|_{10} |n|_{110}
		+ \frac{1}{48} |n|_{10}^3.
\end{multline*}
For slightly larger powers of $2$, the expressions do not become unmanageably large, but the running time of the computation does grow quickly because of the many ways to force $i$ overlaps.
Computing $a_{128}(n)$ took an hour and a half.
After performing these computations, the author was made aware of Everett's work~\cite{Everett 2008}, which raises the possibility of using Everett's expression for $A_n(\beta)$ to compute formulas for $a_{p^\alpha}(n)$ more quickly.
This question has not been investigated, although the results would be interesting to know.

For an integer partition $P$, rewriting $\sum_{\nonoverlappingsubwords{P}{n}} c(v) c(w) \cdots c(z)$ as described yields a multivariate polynomial in $|n|_w$ for various words $w$.
Therefore $a_{p^\alpha}(n)/a_p(n)$ is also a polynomial in $|n|_w$.
For $2 \leq \gamma \leq 2 (\alpha - 1)$, the longest partition in $S_\alpha(\gamma)$ is $\{3, 2, \dots, 2\}$ or $\{2, 2, \dots, 2\}$ and has length $\lfloor \gamma/2 \rfloor$, so the degree of $a_{p^\alpha}(n)/a_p(n)$ is $\alpha - 1$.
Moreover, the longest clusters occurring for a given $\alpha$ have length $\alpha$, so $a_{p^\alpha}(n)/a_p(n)$ depends only on $|n|_w$ for words $w$ of length at most $\alpha$.

\begin{theorem}\label{form}
Let $p$ be a prime, and let $\alpha \geq 1$.
Then $a_{p^\alpha}(n)/a_p(n)$ is a polynomial of degree $\alpha - 1$ in $|n|_w$ for $|w| \leq \alpha$.
\end{theorem}

It would be nice to know more about these polynomials:
How does the number of terms grow?
What can be said about the coefficients?
In particular, why are the coefficients always nonnegative?
Can any sense be made of them as series expansions for $(n+1)/a_p(n)$ if we fix $p$ and let $\alpha \to \infty$?
For example, the coefficient of $|n|_{10}$ in $a_{2^\alpha}(n)/2^{|n|_1}$ for $\alpha = 1, 2, \dots$ takes on the values
\[
	0, \frac{1}{2}, \frac{3}{8}, \frac{5}{12}, \frac{77}{192}, \frac{391}{960}, \frac{259}{640}, \dots,
\]
and a plot of these values suggests that the limit of this sequence exists.

For any $w$ it is known that $|n|_w$ is a $p$-regular sequence in the sense of Allouche and Shallit~\cite[Theorem~6.1]{Allouche--Shallit}.
That is, $|n|_w$ is determined by a finite set of linear recurrences in $|p^e n + i|_w$ along with finitely many initial conditions.
From closure properties of $p$-regular sequences it follows from Theorem~\ref{form} that $a_{p^\alpha}(n)$ is also $p$-regular.
Experimental evidence suggests that the rank of $a_{p^\alpha}(n)$ --- the minimal number of initial conditions required --- is $2 \alpha - 1$.
We leave this as another open problem.

\end{document}